\documentclass{elsart3-1}

\usepackage{amssymb}
\usepackage{amsmath}
\usepackage{enumerate}

\usepackage[T1]{fontenc}
\usepackage[utf8]{inputenc}
\usepackage{lmodern}
\usepackage[english,francais]{babel}

\newtheorem{theorem}{Theorem}[section]

\newtheorem{proposition}[theorem]{Proposition}
\newtheorem{corollary}[theorem]{Corollary}

\renewcommand{\theequation}{\arabic{equation}}
\def\og{\leavevmode\raise.3ex\hbox{$\scriptscriptstyle\langle\!\langle$~}}
\def\fg{\leavevmode\raise.3ex\hbox{~$\!\scriptscriptstyle\,\rangle\!\rangle$}}

\newcommand{\C}{\mathbb{C}}
\def\germ #1 {\mathfrak{#1}}
\def\cal #1 {\mathcal{#1}}

\journal{the Académie des sciences}

\begin{document}
% Vous pouvez mettre dans la prochaine ligne la rubrique choisie
% (si vous la connaissez) - meme deux, format : Rubrique1/Rubrique2
\centerline{}

\begin{frontmatter}

% Titre, auteurs et adresses

% utiliser la commande \thanksref dans \title, \author ou \address
%     pour les notes en bas de page ;

% utiliser la commande \ead pour l'adresse e-mail de chaque auteur
%    (apr� la commande \auteur) ;


 \selectlanguage{francais}

 \title{A generalization of the category $\cal O $ of Bernstein--Gelfand--Gelfand}
 \author{Guillaume Tomasini}
 %%\thanks[label1]{}
 \ead{tomasini@math.u-strasbg.fr}
 \address{IRMA, CNRS et Universit\'e de Strasbourg, \\ 7 rue Ren\'e Descartes, 67084 STRASBOURG CEDEX}

% Vous pouvez mettre a la prochaine ligne les dates
% (de reception et d'acceptation), et le nom du presentateur de votre Note

 \medskip

 %%\begin{center}
 %% {\small Re\c{c}u le *****~; accept\'e apr\`es r\'evision le +++++\\
 %% Pr\'esent\'e par }
 %%\end{center}

 \begin{abstract}
  
 \selectlanguage{francais}
  
	L'\'etude des repr\'esentations irr\'eductibles d'une alg\`ebre de Lie simple d\'efinie sur le corps des nombres complexes a conduit Bernstein, Gelfand et Gelfand a introduire une cat\'egorie qui fournit un cadre naturel pour les modules de plus haut poids. Le but de cette note est de pr\'esenter une construction d'une famille de cat\'egories g\'en\'eralisant celle de Bernstein--Gelfand--Gelfand. Nous d\'ecrivons les modules simples de certaines de ces cat\'egories. Cette classification permet de montrer que ces cat\'egories sont semi--simples.
  
  {\it Pour citer cet article~: }
  \vskip 0.5\baselineskip

  \selectlanguage{english}
  
  \noindent{\bf Abstract}
  \vskip 0.5\baselineskip
  \noindent {\bf }

	In the study of simple modules over a simple complex Lie algebra, Bernstein, Gelfand and Gelfand introduced a category of modules which provides a natural setting for highest weight modules. In this note, we define a family of categories which generalizes the BGG category. We classify the simple modules for some of these categories. As a consequence we show that these categories are semisimple.

  {\it To cite this article:}
 \end{abstract}
\end{frontmatter}

% Maintenant la version abrégée en anglais, si pr�ente

\selectlanguage{english}

%\section*{Abridged English version}
% Texte de la version abrégée� en anglais

\selectlanguage{english}
% texte principale

%%%%%%%%%%%%%%%%%%%%%%%%%%%%%%%%%%%%%%%%%%%%%%%%%%%%%%
\section{Weight modules and Generalized Verma Modules}
%%%%%%%%%%%%%%%%%%%%%%%%%%%%%%%%%%%%%%%%%%%%%%%%%%%%%%

Let $\germ g $ denote a simple Lie algebra over $\C$ and $\cal U (\germ g )$ denote its universal enveloping algebra. Let $\germ h $ be a fixed Cartan subalgebra and denote by $\cal R $ the corresponding set of roots. For $\alpha \in \cal R $ we denote by $H_{\alpha}\in\germ h $ the corresponding co--root and by $\germ g ^{\alpha}$ the rootspace for the root $\alpha$. We will denote by $\Delta$ a set of simple roots. We write $\cal R ^+$ for the corresponding set of positive roots. For a subset $\theta\subset \Delta$, we denote by $\langle \theta \rangle$ the set of all roots which are linear combination of elements of $\theta$ and set $\langle \theta \rangle^{\pm}=\langle \theta \rangle\cap \cal R ^{\pm}$. We consider the following subalgebras of $\germ g $ :
$$\germ l _{\theta}=\germ h \oplus \sum_{\alpha\in \langle \theta\rangle} \: \germ g ^{\alpha}, \quad \germ n ^{\pm}_{\theta}=\sum_{\alpha \in \cal R ^{\pm}-\langle \theta\rangle ^{\pm}} \: \germ g ^{\alpha}.$$

The subalgebra $\germ p _{\theta}=\germ l _{\theta}\oplus \germ n ^+_{\theta}$ is called the standard parabolic subalgebra associated to $\theta$ and $\germ l _{\theta}$ is the standard Levi subalgebra. The later is a reductive algebra. Its semisimple part is denoted $\germ l '_{\theta}$. If $\theta=\emptyset$, then $\germ l _ {\emptyset}=\germ h $ and we simply write $\germ n ^+ $ instead of $\germ n ^+_{\emptyset}$.

We denote by $Mod(\germ g )$ the category of all $\germ g $--modules. We will investigate some full subcategories of $Mod(\germ g )$ for which we will describe the simple modules. A module $M$ is a \emph{weight module} if it is finitely generated, and $\germ h $--diagonalizable in the sense that 
$$M=\oplus_{\lambda \in \germ h ^*}\: M_{\lambda}, \quad M_{\lambda}=\{m\in M \: : \: H\cdot m = \lambda(H)\},$$ with weight spaces $M_{\lambda}$ of finite dimension. We will denote by $\cal M (\germ g , \germ h )$ the full subcategory of $Mod(\germ g )$ consisting of all weight modules. This category has been studied by several authors (e.g. \cite{Fe90}, \cite{Fu}). This category also appears as a particular case of generalised Harisch--Chandra modules (see \cite{PS} for a definition).

The Bernstein--Gelfand--Gelfand category $\cal O $ is the full subcategory of $\cal M (\germ g , \germ h )$ whose objects are $\germ n ^+$--finite (a module $M$ is $\germ n ^+$--finite if for all $m\in M$ the set $\cal U (\germ n ^+)m$ is a finite dimensional vector space). The cateory $\cal O $ was introduced by Bernstein, Gelfand and Gelfand (see \cite{BGG75}). There they desbribe the simple modules and the structure of the category itself, namely they study the projective modules in $\cal O $ and give a certain correspondence for multiplicities (the so called BGG--correspondence). For a review of these results we refer the reader to \cite{Hu08}. In order to discuss the simple modules in $\cal M (\germ g , \germ h )$ we need to recall some well known facts about generalised Verma modules.

A $\germ l _{\theta}$--module can be made into a $\germ p _{\theta}$--module by letting $\germ n ^+_{\theta}$ act trivially. For such a module $N$ we define the \emph{generalised Verma module} (GVM) $V(\theta, N)$ by $V(\theta,N)=\cal U (\germ g )\otimes _{\cal U (\germ p _{\theta})} \: N.$ For any $\germ g $--module $M$ we define the $\germ l _{\theta}$--module $M^{\germ n ^+_{\theta}}:=\{m\in M \: | \: \germ n ^+_{\theta}m=0\}.$ Let us recall the following classical facts about GVMs:

\begin{proposition}\label{GVM}
\begin{enumerate}
\item If $N$ is a simple $\germ l _{\theta}$--module, the module $V(\theta,N)$ admits a unique simple quotient, denoted by $L(\theta,N)$. The module $L(\theta,N)$ is called the simple $\germ g $--module \emph{induced} from $(\germ l _{\theta},N)$.
\item If $M$ is a simple $\germ g $--module such that $M^{\germ n ^+_{\theta}} \not=\{0\}$, then $M\cong L(\theta,M^{\germ n ^+_{\theta}})$.
\end{enumerate}
\end{proposition}
We refer to \cite[proposition $3.8$]{Fe90} for a proof. We refer to \cite{Maz00} for a more detailed discussion of GVM's.

To give the classification of simple weight modules, we need one more ingredient: the so--called \emph{cuspidal} modules.

%%%%%%%%%%%%%%%%%%%%%%%%%%%
\section{Cuspidal modules}
%%%%%%%%%%%%%%%%%%%%%%%%%%%

Let $M$ be a weight module. A root $\alpha \in \cal R $ is said to be \emph{locally nilpotent} (with respect to $M$) if $X_{\alpha}\in \germ g ^{\alpha}$ acts by a locally nilpotent operator on the whole module $M$. It is said to be \emph{cuspidal} if $X_{\alpha}$ acts injectively on the whole module $M$. We denote by $\cal R ^N(M)$ the set of locally nilpotent roots and by $\cal R ^I(M)$ the set of cuspidal roots. We shall simply denote $\cal R ^N $ and $\cal R ^I$ when the module $M$ is clear from the context.

It is known that for a simple weight module $\cal R = \cal R ^N\sqcup \cal R ^I$ (see \cite[lemma $2.3$]{Fe90}). A weight module is called \emph{cuspidal} if $\cal R =\cal R ^I$. Set $\cal R ^N_s=\{\alpha \in \cal R ^N \: : \: -\alpha \in \cal R ^N\}, \quad \cal R ^N_a=\cal R ^N-\cal R ^N_s.$

We define $\cal R ^I_s$ and $\cal R ^I_a$ the same way. Recall the following theorem:

\begin{theorem}[Fernando {\cite[theorem $4.18$]{Fe90}}, Futorny \cite{Fu}]\label{thmFe}%
Let $M$ be a simple weight module. Then $M$ is induced from a cuspidal simple module of some Levi subalgebra of $\germ g $.
\end{theorem}

More precisely there is a set a simple roots $\Delta$ and a subset $\theta\subset \Delta $ such that $\langle \theta \rangle =\cal R ^I_s$ and $\cal R ^+-\langle \theta \rangle ^+=\cal R ^N$. Then $M^{\germ n ^+_{\theta}}$  is a simple cuspidal module for $\germ l _{\theta}$ and $M\cong L(\theta,M^{\germ n ^+_{\theta}})$.

The theorem of Fernando reduces the classification of simple weight $\germ g $--modules to the classification of simple cuspidal weight modules for reductive Lie algebras. By standard arguments this reduces to the classification of simple cuspidal modules for simple Lie algebras. A first step towards this classification is given by the following theorem:

\begin{theorem}[Fernando {\cite[theorem $5.2$]{Fe90}}]
Let $\germ g $ be a simple Lie algebra. If $M$ is a simple cuspidal $\germ g $--module, then $\germ g $ is of type $A$ or $C$
\end{theorem}

The classification of simple cuspidal modules was then completed in two steps. In the first step Britten and Lemire classified all simple cuspidal modules of degree $1$ (see \cite{BL87}) where $deg(M)=\sup_{\lambda \in \germ h ^*} \: \{dim(M_{\lambda})\}.$

Britten and Lemire, and later Benkart, Britten and Lemire, have classified all simple modules of degree $1$ when $\germ g $ is of type $A$ or $C$ (see \cite{BBL97}). These modules will play a key role in our theorem \ref{thmO} below.

Later Mathieu gave the full classification of simple cuspidal modules of finite degree greater than $1$ by introducing the notion of a coherent family (see \cite{Ma00}).

%%%%%%%%%%%%%%%%%%%%%%%%%%%%%%%%%%%%%%%%%%%%
\section{The category $\cal O _{S,\theta}$}
%%%%%%%%%%%%%%%%%%%%%%%%%%%%%%%%%%%%%%%%%%%%

Now we define a family of full subcategories of $\cal M (\germ g ,\germ h )$. Fix a set of simple roots $\Delta$. Fix two subsets $S$ and $\theta$ of $\Delta$ such that $\theta\subset S$. The category $\cal O _{S,\theta}$ is the full subcategory of $\cal M (\germ g , \germ h )$ whose objects are the modules satisfying the following conditions :
\begin{enumerate}[$(\mathcal{O}1)$]\label{defO}
\item As a $\germ l _{\theta}$--module, $M$ is a direct sum of simple highest weight modules.\label{O1}
\item For $\alpha \in \langle S-\theta\rangle $, the element $X_{\alpha}\in \germ g ^{\alpha}$ is cuspidal for $M$.\label{O2}
\item The module $M$ is $\germ n ^+_S$--finite.\label{O3}
\end{enumerate}

Notice than in the case $S=\theta=\emptyset$, we recover the definition of the category $\cal O $. Other generalisations of category $\cal O $ can be recovered for particular choices of $(\theta,S)$ (see \cite{RC80}, \cite{Maz00}). We sometimes write $\cal O _{S,\theta}(\germ g )$ to emphasize the Lie algebra. Now we list some easy properties of $\cal O _{S,\theta}$ which are generalisations of analogous properties of $\cal O $.

\begin{proposition}\label{propO}
Let $M,N \in \cal O _{S,\theta}$.
\begin{enumerate}
\item Then $M\oplus N$ is in $\cal O _{S,\theta}$. Every submodule and every quotient of $M$ is again in $\cal O _{S,\theta}$.
\item The category $\cal O _{S,\theta}$ is abelian, and every module in $\cal O _{S,\theta}$ is noetherian and artinian.
\item The category $\cal O _{S,\theta}$ is $\cal Z (\germ g )$--finite (where $\cal Z (\germ g )$ is the center of $\cal U (\germ g )$).
\item The simple modules in $\cal O _{S,\theta}$ are modules of the form $L(S,N)$ where $N$ is a simple module in $\cal O _{S,\theta}(\germ l _S)$.
\end{enumerate}
\end{proposition}

Being artinian and noetherian, every module $M$ in $\cal O _{S,\theta}$ admits a finite Jordan--H\"older series whose quotients are of the form $L(S,N)$ (see \cite{Ja80}). This allows us to define the multipicity of $L(S,N)$ in any Jordan--H\"older series of $M$. From proposition \ref{propO}$(iv)$, the description of the simple modules in $\cal O _{S,\theta}(\germ g )$  first requires to consider the case where $S=\Delta$.

In the sequel we assume $S=\Delta$. We shall only consider a simple Lie algebra $\germ g $ and restrict ourselves to the case $\theta \not=\emptyset$ and $\theta\not=\Delta$. We label the simple roots as in \cite{Bo68}. First, from theorem \ref{thmFe}, every simple module in $\cal O _{\Delta,\theta}$ is of the form $L(\Delta-\theta,C)$ for some simple cuspidal $\germ l _{\Delta-\theta}$--module $C$. Such a module is defined by a simple cuspidal $\germ l '_{\Delta-\theta}$--module $C'$ and a scalar action of the center $\germ z $ of $\germ l _{\Delta-\theta}$. We have to find those $C$ for which $L(\Delta-\theta,C)$ satisfies condition $(\mathcal{O}1)$ of the definition of the category $\cal O _{\Delta,\theta}$.

\begin{theorem}\label{thmO}
Assume $(\germ g ,\Delta-\theta)$ is not one of the following: $(B_n,\{e_1\})$, $(D_n,\{e_1\})$, $(D_n,\{e_{n-1}\})$, $(D_n,\{e_n\})$, $(E_6,\{e_1\})$, $(E_6,\{e_6\})$, $(E_7,\{e_7\})$. Then we have
\begin{enumerate}
\item If $L(\Delta-\theta,C)$ is a simple module in $\cal O _{\Delta,\theta}(\germ g )$, then $\germ g $ is of type $A_n$ or is of type $C_n$, $\Delta-\theta$ is a connected subset of $\Delta$ and $C$ is a simple cuspidal module of degree $1$.
\item If $\germ g =\germ sl _{n+1}$ and $\Delta-\theta$ is any connected subset of $\Delta$ other than $\{e_1\}$ and $\{e_n\}$, then the simple modules in $\cal O _{\Delta,\theta}(\germ g )$ are modules of degree $1$. Conversely, any infinite dimensional simple $\germ sl _n$--module of degree $1$ is a simple object in some category $\cal O _{\Delta,\theta}(\germ sl _n)$.
\item If $\germ g =\germ sp _n$, then there are simple modules in $\cal O _{\Delta,\theta}(\germ g )$ if and only if $\Delta-\theta\supset \{e_n\}$. In this case, all the simple modules in $\cal O _{\Delta,\theta}(\germ g )$ are modules of degree $1$. Conversely, any infinite dimensional simple $\germ sp _n$--module of degree $1$ is a simple object in some category $\cal O _{\Delta,\theta}(\germ sp _n)$.
\end{enumerate}
\end{theorem}

From now on, we assume that $\germ g $ is of type $A$ or $C$ and that $\theta$ is neither empty nor equal to $\Delta$. When $\germ g = A_n$ we also assume that $\Delta-\theta$ is different from $\{e_1\}$ and $\{e_n\}$. When the category $\cal O _{\Delta,\theta}$ is non empty, we want to describe its structure. First we note the following fact: if $F$ is a finite dimensional $\germ g $--module of dimension greater than $1$ and $M\in \cal O _{\Delta,\theta}$ is simple, then $F\otimes M$ does not belong to the category $\cal O _{\Delta,\theta}$ in general. The proof of this fact requires some results of Britten and Lemire \cite[section $3$]{BL01}. Therefore the category $\cal O _{\Delta,\theta}$ has a structure very different from the case of the category $\cal O $. In fact, we prove the following

\begin{theorem}\label{thmext}
Let $M$ and $N$ be two simple modules in $\cal O _{\Delta,\theta}$. Then $Ext^1_{\cal O _{\Delta,\theta}}(M,N)=0$.
\end{theorem}

As a corollary we get:

\begin{corollary}
The category $\cal O _{\Delta,\theta}$ is semisimple.
\end{corollary}

Remark that the result holds trivially when $\theta=\Delta$ (this is in fact a part of the definition of the category). Note that when $\theta=\emptyset$ and $\germ g $ is of type $C$, Britten, Khomenko, Lemire, Mazorchuk have proved the semisimplicity of $\cal O _{\Delta,\theta}$ in \cite{BKLM}. The result does not hold when $\theta=\emptyset$ and $\germ g $ is of type $A$. In \cite{GS}, Grantcharov and Serganova have described the indecomposable modules in this case. The proof of theorem \ref{thmext} uses the example $3.3$ of \cite{GS}.

Detailed proofs will be published elsewhere.

%%%%%%%%%%%%%%%%%%%%%%%%%%%%%%%%%%%%%%%%%%%%%%%%%%%%%%%


\begin{thebibliography}{00}
% Essayez �utiliser le systeme 'bibitem',
%    avec les references en ordre alphabetique.
% Texte


 \bibitem{BBL97}
 G. Benkart, D. Britten, F. Lemire,
 Modules with bounded weight multiplicities for simple Lie algebras,
 Math. Z.
 (2) 225 (1997) 333-353.
 
 \bibitem{BGG75} 
 I. N. Bern{\v{s}}te{\u\i}n, I. M. Gel$'$fand, S. I. Gel$'$fand,
 Differential operators on the base affine space and a study of {${\germ g }$}-modules,
 in Lie groups and their representations (Proc. Summer School, Bolyai J\'anos Math. Soc., Budapest, 1971)
 (1975) 21-64.
 
 \bibitem {Bo68}
 N. Bourbaki,
 Groupes et {A}lg\`ebres de {L}ie, Chap. IV--VI,
 Hermann, Paris, 1968.

 \bibitem{BKLM}
 D. Britten, O. Khomenko, F. Lemire, V. Mazorchuk,
 Complete reducibility of torsion free {$C_n$}-modules of finite degree,
 J. Algebra
 276 (2004) 129-142.

\bibitem{BL87}
 D. Britten, F. Lemire,
 A classification of simple {L}ie modules having a {$1$}-dimensional weight space,
 Trans. Amer. Math. Soc..
 (2) 299 (1987) 683-697.

\bibitem{BL01}
  D. Britten, F. Lemire,
  Tensor product realizations of simple torsion free modules,
  Canad. J. Math.,
  (2) 53 (2001) 225--243.

\bibitem{Fe90}
 S. L. Fernando,
 Lie algebra modules with finite-dimensional weight spaces. {I},
 Trans. Amer. Math. Soc.
 (2) 322 (1990) 757-781.

\bibitem{Fu}
 V. Futorny,
 The weight representations of semisimple finite dimensional Lie algebras,
 PhD Thesis,
 Kiev University,
 1987.
 
\bibitem{GS}
 D. Grantcharov, V. Serganova,
 Cuspidal representations of {$\mathfrak{sl}(n+1)$},
 arXiv:0710.2682v1.

\bibitem{Hu08}
 J. E. Humphreys,
 Representations of semisimple {L}ie algebras in the {BGG} category {$\mathcal{O}$},
 vol. 94 of Graduate Studies in Mathematics,
 American Mathematical Society, Providence, RI, 2008.

\bibitem{Ja80}
 N. Jacobson,
 Basic Algebra $II$,
 W.H. Freeman, 1980.

\bibitem{Ma00}
 O. Mathieu,
 Classification of irreducible weight modules,
 Ann. Inst. Fourier (Grenoble),
 (2) 50 (2000) 537-592.

\bibitem{Maz00}
 V. Mazorchuk,
 Generalized {V}erma modules,
 vol. 8 of Mathematical Studies Monograph Series,
 VNTL Publishers, L$'$viv, 2000.

\bibitem{PS}
 I. Penkov, V. Serganova,
 Generalized {H}arish-{C}handra modules,
 Mosc. Math. J.,
 (2) (2002) 753-767.

\bibitem{RC80}
  A. Rocha-Caridi,
  Splitting criteria for $\germ g $-modules induced from a parabolic and the {B}er\v nste\u\i n-{G}el$'$fand-{G}el$'$fand resolution of a finite-dimensional, irreducible {${\germ g }$}-module,
  Trans. Amer. Math. Soc.,
  (2) 262 (1980) 335--366.

%\bibitem{To09}
% G. Tomasini,
% A {H}owe--type correspondence for the dual pair $(\germ sl _2,\germ sl _n)$ in $\germ sl _{2n}$.
% arXiv:0911.5210

\end{thebibliography}
\end{document}